\documentclass[12pt]{amsart}

\usepackage{mathrsfs}

\usepackage{amsmath,amsthm,amssymb}
\font\teneufm=eufm10
\font\seveneufm=eufm7
\font\fiveeufm=eufm5
\newfam\frakturfam

\textfont\frakturfam=\teneufm
\scriptfont\frakturfam=\seveneufm
\scriptscriptfont\frakturfam=\fiveeufm

\newtheorem{pr}{Proposition}

\newtheorem{lemma}{Lemma}
\newtheorem{theorem}{Theorem}
\newtheorem{corol}{Corollary}

\def\bee{\begin{eqnarray}}
\def\bes{\begin{eqnarray*}}
\def\eee{\end{eqnarray}}
\def\ees{\end{eqnarray*}}

\def\a{\alpha}
\def\b{\beta}

\def\Proof{{\sl Proof.}\ }

\newcommand{\Aut}{\mbox{Aut}}

\title{Automorphisms of simple quotients of the Poisson and universal enveloping algebras of $\mathrm{sl}_2$}

\begin{document}
\date{}

\maketitle

\bigskip
\begin{center}

{\bf Altyngul Naurazbekova}
\footnote{Department of Mathematics, L.N. Gumilyov Eurasian National University, Nur-Sultan, Kazakhstan,  e-mail: {\em altyngul.82@mail.ru}},
{\bf Ualbai Umirbaev}\footnote{Department of Mathematics,
 Wayne State University,
Detroit, MI 48202, USA; Department of Mathematics, 
Al-Farabi Kazakh National University, Almaty, 050040, Kazakhstan; 
and Institute of Mathematics and Mathematical Modeling, Almaty, 050010, Kazakhstan,
e-mail: {\em umirbaev@wayne.edu}}
\end{center}

\begin{abstract} Let $P(\mathrm{sl}_2(K))$ be the Poisson enveloping algebra of the Lie algebra $\mathrm{sl}_2(K)$ over an algebraically closed field $K$ of characteristic zero. The quotient algebras $ $ $P(\mathrm{sl}_2(K))/(C_P-\lambda)$, where $C_P$ is the standard Casimir element of $\mathrm{sl}_2(K)$ in $P(\mathrm{sl}_2(K))$ and  $0\neq \lambda\in K$, are proven to be simple in \cite{UZh}. Using a result by L. Makar-Limanov \cite{ML90}, we describe generators of the automorphism group of $P(\mathrm{sl}_2(K))/(C_P-\lambda)$ and represent this group as an amalgamated product of its subgroups. Moreover, using similar results by J. Dixmier \cite{Dixmier73} and O. Fleury \cite{Fleury} for the quotient algebras  $U(\mathrm{sl}_2(K))/(C_U-\lambda)$, where $C_U$ is the standard Casimir element of $\mathrm{sl}_2(K)$ in the universal enveloping algebra $U(\mathrm{sl}_2(K))$, we prove that the automorphism groups of $P(\mathrm{sl}_2(K))/(C_P-\lambda)$ and $U(\mathrm{sl}_2(K))/(C_U-\lambda)$ are isomorphic. 
\end{abstract}



\section{Introduction}
\hspace*{\parindent}

It is well known that all automorphisms of the polynomial algebra $K[x,y]$ in two variables $x,y$ over a field $K$ are tame \cite{Jung, Kulk}. Moreover, the automorphism group $ \mathrm{Aut}\,K[x,y]$ of this algebra admits an amalgamated free product structure \cite{Kulk, Shafar}, i.e.,
\bes
\mathrm{Aut}\,K[x,y]=A\ast_{C}B, 
\ees
where $A$ is the affine automorphism subgroup, $B$ is the triangular automorphism subgroup, and  $C=A\cap B$. 

Similar results hold for free associative algebras \cite{Czer, ML70},  free Poisson algebras (in characteristic zero) \cite{MLTU}, and free right-symmetric algebras of rank two \cite{ANK,KMLU}.  Moreover, the automorphism groups of polynomial algebras, free associative algebras, and free Poisson algebras in two variables are isomorphic. 

The automorphism  groups of commutative and associative algebras generated by three elements are much more complicated. 
The well-known Nagata automorphism (see \cite{Nagata})
\bes
\sigma=(x+2y(zx-y^2)+z(zx-y^2)^2, y+z(zx-y^2),z)
\ees
of the polynomial algebra $K[x,y,z]$ over a field $K$ of characteristic $0$ is proven to be non-tame \cite{Umi25}. Free associative algebras  in three variables over a
field of characteristic $0$ also have non-tame automorphisms \cite{Umi33}.
The Nagata
automorphism gives an example of a wild
automorphism of free Poisson algebras in three variables.

 In 1964 P. Cohn  proved \cite{Cohn} that all automorphisms of finitely generated free
Lie algebras over a field are tame.
 Later this result was extended to free algebras of Nielsen-Schreier varieties \cite{Lewin}. Recall that
a variety of universal algebras is called Nielsen-Schreier, if any subalgebra of a free algebra of this variety is free,
i.e., an analog of the classical Nielsen-Schreier theorem is true.
 The varieties of all non-associative algebras \cite{Kurosh},
 commutative and anti-commutative algebras \cite{Shirshov54}, Lie algebras
 \cite{Shirshov53,Witt}, and Lie superalgebras  \cite{Mikhalev,Stern} over a field are Nielsen-Schreier. 

The automorphism groups  of free non-associative algebras and free commutative algebras of rank two admit an amalgamated free product structure \cite{ANK}. The groups of automorphisms of free Lie algebras and free anti-commutative algebras of rank three also admit an amalgamated free product structure \cite{ANU20}.

The study of relations between the automorphism groups of Poisson algebras and their deformation-quantizations is motivated by the {\em Belov-Kanel -- Kontsevich Conjecture} \cite{BKK2}, which asserts that the automorphism group of the Weyl algebra $A_n$ of index $n\geq 1$ over a field $K$ is isomorphic to the group of automorphisms of the symplectic Poisson algebra $P_n$, i.e.,
\bes
\Aut\,A_n \simeq \Aut\,P_n.
\ees
Recall that the Weyl algebra $A_n$ of index $n\geq 1$,  is the associative algebra over a field $K$ with generators $X_1,\ldots,X_n,Y_1,\ldots,Y_n$ and defining relations
\bes
[Y_i,X_j]=\delta_{ij}, \ \ [X_i,X_j]=0, \ \
[Y_i,Y_j]=0,
\ees
where $\delta_{ij}$ is the Kronecker symbol
and $1\leq i,j\leq n$. The symplectic Poisson algebra $P_n$ of index $n\geq 1$ is the polynomial algebra in the variables
$x_1,\ldots,x_n,y_1,\ldots,y_n$ endowed with the Poisson bracket
defined by
\bes
\{y_i,x_j\}=\delta_{ij}, \ \ \{x_i,x_j\}=0, \ \
\{y_i,y_j\}=0,
\ees
where $1\leq i,j\leq n$.

The structure of the automorphism group of the Weyl algebra $A_1$ is described in  \cite{Alev86,Dixmier68,ML84}. These results easily imply that the groups of automorphisms of the symplectic Poisson algebra $P_1$ and the Weyl algebra $A_1$ are isomorphic, i.e., 
the Belov-Kanel -- Kontsevich Conjecture is true for $n=1$. In 2005 Belov-Kanel and Kontsevich \cite{BKK2} proved that the groups of tame automorphisms of the symplectic Poisson algebra $P_n$ and the Weyl algebra $A_n$ are isomorphic.

Let  $L$  be an arbitrary Lie algebra over a field $K$ of characteristic zero. Denote by $U(L)$ the universal enveloping algebra of $L$ and by $P(L)$ the Poisson enveloping algebra of $L$ \cite{UZh} (or the symmetric Poisson algebra of $L$). 
Recall that the well known symmetrization map \cite[p.77]{Dixmier}
\bes
 S: P(L)\rightarrow U(L)
\ees  is an isomorphism of $L$-modules. Moreover, $U(L)$ is the most well-known natural deformation-quantization of $P(L)$.

 In 1973 J. Dixmier  \cite{Dixmier73} studied the quotients $U_{\lambda}=U(\mathrm{sl}_2(\mathbb{C}))/(C_U-\lambda)$ of the universal enveloping algebra $U(\mathrm{sl}_2(\mathbb{C}))$ of the three dimensional simple Lie algebra $\mathrm{sl}_2(\mathbb{C})$ over the field of complex numbers $\mathbb{C}$, where $C_U$ is the standard Casimir element and $0\neq \lambda \in \mathbb{C}$. The structure of these algebras depend on $\lambda$ and $U_{\lambda}$ is simple if $\lambda\neq n^2+2n$ for any natural $n$. If $\lambda=n^2+2n$, then $U_{\lambda}$ has a unique non-trivial ideal of finite codimension $(n+1)^2$. In the same paper Dixmier described generators of the group of automorphisms of the algebra $U_{\lambda}$ and defined the group of tame automorphisms of $U(\mathrm{sl}_2(\mathbb{C}))$. He also formulated a question on the existence of wild automorphisms of 
$U(\mathrm{sl}_2(\mathbb{C}))$. In 1976 A. Joseph \cite{Joseph}  gave an example of a wild automorphism of $U(\mathrm{sl}_2(\mathbb{C}))$. In 1998 O. Fleury \cite{Fleury} represented the automorphism group of the quotient algebra $U_{\lambda}$ as an amalgamated product of its subgroups. Using this, she also proved that every finite subgroup of $\mathrm{Aut}\,U(\mathrm{sl}_2(\mathbb{C}))$ is isomorphic to a subgroup of $\mathrm{Aut}\,\mathrm{sl}_2(\mathbb{C})$.

Let $P(\mathrm{sl}_2(K))$ be the Poisson enveloping algebra of the Lie algebra $\mathrm{sl}_2(K)$ over an algebraically closed field $K$ of characteristic zero. The quotient algebras $P_{\lambda}=P(\mathrm{sl}_2(K))/(C_P-\lambda)$, where $C_P$ is the standard Casimir element of $\mathrm{sl}_2(K)$ in $P(\mathrm{sl}_2(K))$ and  $0\neq \lambda\in K$, are proven to be simple in \cite{UZh}. Notice that the Casimir elements $C_P$ of $P(\mathrm{sl}_2(K))$ and $C_U$ of $U(\mathrm{sl}_2(K))$  correspond to each other under the symmetrization map. 

In 1990 L. Makar-Limanov \cite{ML90} described generators of the automorphism group of the algebraic surface defined by $xy=f(z)$ over an algebraically closed field. Using this result, we describe generators of the automorphism group  of the Poisson quotient algebra $P(\mathrm{sl}_2(K))/(C_P-\lambda)$, where $\lambda \in K$, and represent this group as an amalgamated product of its subgroups. Then, using the above described results by J. Dixmier \cite{Dixmier73} and O. Fleury \cite{Fleury}, we prove that the automorphism groups of $P(\mathrm{sl}_2(K))/(C_P-\lambda)$ and $U(\mathrm{sl}_2(K))/(C_U-\lambda)$ are isomorphic. 

Unfortunately, the question on the isomorphism of the automorphism groups of $P(\mathrm{sl}_2(K))$ and $U(\mathrm{sl}_2(K))$ remains open. 

The paper is organized as follows. In Section 2, we describe generators of the automorphism group of the Poisson quotient algebra $P(\mathrm{sl}_2(K))/(C_P-\lambda)$, where $\lambda \in K$. In Section 3,  we prove that the automorphism group of $P(\mathrm{sl}_2(K))/(C_P-\lambda)$  admits an amalgamated free product structure.  In Section 4, we show that the automorphism groups of $P(\mathrm{sl}_2(K))/(C_P-\lambda)$ and $U(\mathrm{sl}_2(K))/(C_U-\lambda)$ are isomorphic.

\section{ Generators of the automorphism group of $P_{\lambda}$}
\hspace*{\parindent}

A vector space $P$ over a field $K$ endowed with two bilinear operations $x\cdot y$  (a multiplication) and $\{x,y\}$ (a Poisson bracket) is called a \emph{Poisson algebra} if $P$ is a commutative associative algebra under $x\cdot y$, $P$ is a Lie algebra under  $\{x,y\}$, and $P$ satisfies the following identity (the Leibniz identity):
$$ \{x,y\cdot z\}=y\cdot \{x,z\}+ \{x,y\}\cdot z.$$

Let  $L$  be an arbitrary Lie algebra with Lie bracket $[\,,]$ over a field $K$ and let
$e_1,e_2\ldots$ be a linear basis of $L$. Then there exists a unique bracket $\{\,,\}$ on the
 polynomial algebra
$K[e_1, e_2,\ldots]$  defined by $\{e_i,e_j\}=[e_i,e_j]$ for all $i,j$ and satisfying the Leibniz identity.  
With this bracket   
\bes
P(L)=\langle K[e_1,e_2,\ldots],\cdot, \{\,,\}\rangle
\ees
 becomes a Poisson algebra. This Poisson algebra $P(L)$ is called the {\em Poisson enveloping algebra} of $L$ \cite{UZh}. Note that the bracket $\{,\}$ of the algebra $P(L)$ depends on the structure of $L$ but does not depend on a chosen basis.

Consider the three dimensional simple Lie algebra  $\mathrm{sl}_{2}(K)$ over an algebraically closed field $K$ of characteristic zero. 
Let 
\bes
 E=\left[\begin{array}{cc}
 0  &  1\\
 0  &  0\\
\end{array}\right], 
H=\left[\begin{array}{cc}
 1  &  0\\
 0  &  -1\\
\end{array}\right], 
F=\left[\begin{array}{cc}
 0  &  0\\
 1  &  0\\
\end{array}\right] 
\ees
be  the standard basis of the Lie algebra $\mathrm{sl}_{2}(K)$. We have
$$[E,F]=H,  [H,E]=2E, [H,F]=-2F.$$

The center of the Poisson enveloping algebra $P(\mathrm{sl}_2(K))$ of the Lie algebra $\mathrm{sl}_2(K)$ is equal to $K[C_P]$, where 
$$C_P=4EF+ H^{2}$$
 is the Casimir element of $\mathrm{sl}_2(K)$ in $P(\mathrm{sl}_2(K))$. This is an easy corollary of the fact that the center of the universal enveloping algebra $U(\mathrm{sl}_{2}(K))$ of $\mathrm{sl}_{2}(K)$ is generated by the 
standard Casimir element  
$$C_U=4FE+H^{2}+2H=4EF+H^{2}-2H=2EF+2FE+H^{2}$$
 of the universal enveloping algebra $U(\mathrm{sl}_{2}(K))$ of $\mathrm{sl}_{2}(K)$ and the image of $C_P$ under the symmetrization map is $C_U$ \cite{Dixmier}.

For any $\lambda\in K$ let  $P_{\lambda}=P(\mathrm{sl}_{2}(K))/(C_P-\lambda)$ be the quotient algebra of the algebra $P(\mathrm{sl}_{2}(K)$ by the principal ideal $(C_P-\lambda)$. 
 Denote by  $e,h,f$  the images of $E,H,F$ in $P_{\lambda}$, respectively. Then we have
\begin{equation}\label{f1}
\{e,f\}=h, \{h,e\}=2e, \{h,f\}=-2f, 4ef+h^{2}=\lambda. 
\end{equation}

Notice that $P_{\lambda}$, as an associative and commutative algebra, is generated by   $e,h,f$  and defined by one relation 
\begin{equation}\label{f2}
ef=-\frac{1}{4}h^{2}+\frac{1}{4}\lambda.
\end{equation}
Consequently, the set of words of the form 
\bee\label{f3}
 f^{m}h^{n}, \ \ h^{n}e^{r}, \ \  m,n,r\geq 0
\eee
is a linear basis of $P_{\lambda}$.

Let  $G=\mathrm{Aut}(P_{\lambda})$ be the automorphism group of the algebra $P_{\lambda}$. Denote by $\varphi=(f_1,f_2,f_3)$  the automorphism of the algebra $P_{\lambda}$ such that  $\varphi(e)=f_1, \varphi(h)=f_2, \varphi(f)=f_3$. If $\theta=(f_1,f_2,f_3)$ and $\varphi=(g_1,g_2,g_3)$, then
 the product in $\mathrm{Aut}(P_{\lambda})$ is defined by 
$$\theta\circ \varphi=(g_1(f_1,f_2,f_3),g_2(f_1,f_2,f_3),g_3(f_1,f_2,f_3)).$$

First we describe the linear automorphisms of $P_{\lambda}$, i.e., the automorphisms of $\mathrm{sl}_{2}(K)$.
Let $A=\mathrm{Aut}(\mathrm{sl}_{2}(K))$ be the group of all automorphisms of the Lie algebra $\mathrm{sl}_{2}(K)$. Every automorphism of $\mathrm{sl}_{2}(K)$ gives a unique automorphism of $\mathrm{Aut}(P_{\lambda})$. Further we identify $A$ with the corresponding subgroup of $G$. 

It is well known \cite[p. 306]{Jacob} that every automorphism of the algebra $\mathrm{sl}_{2}(K)$ is inner, i.e., every automorphism coincides with 
\bes
\widehat{T} : \mathrm{sl}_{2}(K)\rightarrow \mathrm{sl}_{2}(K),  X\mapsto T^{-1}XT
\ees 
for some matrix $T\in \mathrm{GL}_2(K)$. Consequently, $A\simeq PSL_{2}(K)=GL_{2}(K)/ \{\a I | \a\in K^*\}=SL_{2}(K)/ \{I,-I\}$, where $I$ is the identity matrix. 

Let $C$ be the set of all automorphisms $\widehat{T}$, where $T$ runs over the set of matrices of the form 
\bee\label{f4}
T=\left[\begin{array}{cc}
 1  &  0\\
 \a  &  \b\\
\end{array}\right]\in GL_{2}(K). 
\eee 
Identifying the matrices $E,H,F$ with their images $e,h,f$, we can write  
\bes
\widehat{T}= (\b e+\a h- \a^2/\b f, h-2\a/\b f,  1/\b f). 
\ees
The elements of $C$ will be called {\em linear triangular} automorphisms. Denote by $H$ the subset of $C$ with $\a=0$ and by $C_1$ the subset of $C$ with $\b=1$. This means that $H$ consists of all automorphisms of the form 
\bes
H_\b=(\b e, h,1/\b f) 
\ees
and $C_1$ consists of all automorphisms of the form 
\bes
(e+\a h- \a^2 f, h-2\a f,  f). 
\ees
The automorphisms $H_\b$ are called {\em hyperbolic rotations} \cite{ML90} and the subgroup $H$ will be called the group of {\em hyperbolic rotations}.  
It is easy to check that
$$C\simeq C_1 \rtimes H.$$

\begin{lemma} \label{l1} The system of elements
$$A_0=\{\tau_{\a}=(f,-h+2\a f,e+\a h -\a^{2}f), \mathrm{id}=(e,h,f) \}$$
is a system of representatives of the left cosets of  $C$ in $A$. 
\end{lemma}
\Proof Let 
\bes
T=\left[\begin{array}{cc}
 a  &  b\\
 c  &  d\\
\end{array}\right]\in GL_{2}(K). 
\ees
If $b=0$, then $\widehat{T}\in C$ and $\mathrm{id}$ is a representative of the left coset $\widehat{T}\circ C=\mathrm{id}\circ C=C$. If $b\neq 0$, then we can assume $b=1$ in $PSL_{2}(K)$, i.e., 
\bes
T=\left[\begin{array}{cc}
 a  &  1\\
 c  &  d\\
\end{array}\right]. 
\ees
This matrix can be easily represented as $T=PQ_{\a}$, where 
\bes
P=\left[\begin{array}{cc}
 a_1  &  0\\
 c_1 &  d_1\\
\end{array}\right], \ \ 
Q_{\a}=\left[\begin{array}{cc}
 \a  &  1\\
 1  &  0\\
\end{array}\right]. 
\ees
Then $\widehat{T}=\widehat{Q_{\a}}\circ \widehat{P}\in \widehat{Q_{\a}}\circ C$. Notice that $\widehat{Q_{\a}}=\tau_{\a}$. 

Suppose that $\widehat{Q_{\a_1}}\circ C=\widehat{Q_{\a_2}}\circ C$. This means that 
\bes
\widehat{Q_{\a_1}}^{-1}\circ \widehat{Q_{\a_2}}=\widehat{Q_{\a_2}Q_{\a_1}^{-1}}\in C. 
\ees
We have 
\bes
Q_{\a_2}Q_{\a_1}^{-1}=\left[\begin{array}{cc}
 \a_2  &  1\\
 1  &  0\\
\end{array}\right]
\left[\begin{array}{cc}
 0  &  1\\
 1  &  -\a_1\\
\end{array}\right]=
\left[\begin{array}{cc}
 1  &  \a_2-\a_1\\
 0 &  1\\
\end{array}\right].
\ees
Obviously, $\widehat{Q_{\a_2}Q_{\a_1}^{-1}}\in C$ if and only if $\a_1=\a_2$.  $\Box$

Generators of the automorphism group of the quotient algebra $R=K[x,y,z]/(xy-P(z))$ of the polynomial algebra $K[x,y,z]$ over an algebraically closed field $K$, where $P(z)\in K[z]$, were given by L. Makar-Limanov \cite{ML90}.

\begin{theorem} \label{th_ML} \cite[p. 252]{ML90} The group $\mathrm{Aut}(R)$ is generated by the following automorphisms:

$(a)$ Hyperbolic rotations $H(x)=\nu x, H(y)=\nu^{-1}y, H(z)= z $, $\nu\in K^{\ast}$;

$(b)$ Involution $V(x)=y, V(y)=x, V(z)=z $;

$(c)$ Triangular automorphisms $\Delta(x)=x, \Delta(y)=y+[P(z+xg(x))-P(z)]x^{-1}, \Delta (z)=z+xg(x) $, $g(x)\in K[x]$;

$(d)$ If $P(z)=c(z+a)^{d}$, then rescalings $R(x)=x, R(y)=\nu^{d}y, R(z)=\nu z +(\nu-1)a$, $\nu\in K^{\ast}$ should be added;

$(e)$ If $P(z)=(z+a)^{i}Q((z+a)^{n})$ and $\mu\in K$ such that $\mu^{n}=1$, then the symmetry  $S(x)=x, S(y)=\mu^{i}y, S(z)=\mu z+(\mu-1)a$ should be added;

$(f)$ If $char K=p>0$ and $P(z)=Q(z^{p}-a^{p-1}z)$, then the translation $T(x)=x, T(y)=y, T(z)=z-a$ should be added.

\end{theorem}

\begin{pr} \label{p1}  The group $G=\mathrm{Aut}(P_{\lambda})$ is generated by 
\bes
\tau(e)=f, \tau(h)=-h, \tau(f)=e 
\ees
and the automorphisms of the form 
\begin{equation}\label{5}
\begin{split}
\Delta_g(e)=e-g(f)h-fg^{2}(f),\\
\Delta_g(h)=h+2fg(f),\\
\Delta_g(f)=f , 
\end{split}
\end{equation}
where $g(x)\in K[x]$. 
\end{pr}
\Proof Notice that $P_{\lambda}$ as an associative and commutative algebra  is represented as $R=K[f,e,h]/(fe-P(h))$, where $P(h)=-\frac{1}{4}h^{2}+\frac{\lambda}{4}$. Let $G'=\mathrm{Aut}(R)$ be the group of all automorphisms of the associative and commutative algebra $R$. Consider the list of generators of the group $G'$ given by Theorem \ref{th_ML} under the correspondence $(x,y,z)\rightarrow (f,e,h)$: 

$(a)$ Hyperbolic rotations $H_\nu(e)=\nu e, H_\nu(h)= h, H_\nu(f)=\nu^{-1} f $, $\nu\in K^{\ast}$. 

$(b)$ Involution $V(e)=f,  V(h)=h, V(f)=e$.

$(c)$ Triangular automorphisms 
\bes
\Delta(e)=e+[P(h+fg(f))-P(h)]f^{-1}, \Delta (h)=h+fg(f), \Delta(f)=f,
\ees
where $g(x)\in K[x]$. Notice, that 
\bes
[P(h+fg(f))-P(h)]f^{-1}=-1/2  g(f)h-1/4 f g^{2}(f).
\ees 
Replacing $1/2 g(f)$ by $g(f)$, we get automorphisms $\Delta_g$ of the form \eqref{5}.

$(d)$ Notice that $P(h)$ can be represented in the form $P(h)=c(h+a)^{d}$ only if $\lambda=0$. If $\lambda=0$, then $P(h)=-\frac{1}{4}h^{2}$ and rescalings $ R_{\nu}(e)=\nu^{2}e, R_{\nu}(h)=\nu h , R_{\nu}(f)=f$, where $\nu\in K^{\ast}$, are automorphisms of $R$.

$(e)$ Consider representations of $P(h)$ in the form $P(h)=(h+a)^{i}Q((h+a)^{n})$. If $n=1$ Theorem \ref{th_ML}(e) gives only the identity automorphism. If $n=2$, then, obviously,  $i=0$, $a=0$, and  $\mu=\pm 1$. This gives the automorphism $S(e)=e,  S(h)=-h, S(f)=f$.

However, not all automorphisms from this list are automorphisms of the Poisson algebra $P_{\lambda}$. More exactly, $V,S, R_\nu\notin G$ and $\Delta_g, H_\nu\in G$.

Notice that 
\bee\label{f6}
V^2=S^2=\mathrm{id}
\eee
and that 
\begin{equation}\label{6}
\tau=V\circ S=S \circ V=(f,-h,e)
\end{equation}
is an automorphism of the Poisson algebra $P_{\lambda}$.

Denote by $G_1$ the subgroup of $G$ generated by $\tau$ and all automorphisms of the form \eqref{5}.

We also have 
\bes H_\nu=(\nu e,h,\nu^{-1} f)=\tau\circ (e-i\sqrt{\nu}h+\nu f,h+2i\sqrt{\nu}f,f)\circ \tau\circ\\
\circ (e-i\sqrt{\nu}^{-1} h+\nu^{-1} f,h+2i\sqrt{\nu}^{-1} f,f)\circ \tau\circ (e-i\sqrt{\nu} h+\nu f,h+2i\sqrt{\nu} f,f), 
\ees 
where $i=\sqrt{-1}$, i.e., $H_\nu\in G_1$.  

It is easy to check that
\bee\label{f7}
\Delta_{g_1}\circ \Delta_{g_2}=\Delta_{g_1+g_2}, \ \ \Delta_g^{-1}=\Delta_{-g}. 
\eee

Direct calculations give that 
\bes
S\circ \Delta_g \circ S=(e+g(f)h-fg^{2}(f),-h+2fg(f),f) \circ (e,-h,f)\\
=(e+g(f)h-fg^{2}(f),h-2fg(f),f) =\Delta_{-g}
\ees
and
\bes
 \tau \circ \Delta_{-g} \circ \tau=(f-g(e)h-eg^{2}(e),-h-2eg(e),e)\circ (f,-h,e)\\
=(e,h+2eg(e),f-g(e)h-eg^{2}(e) ) \\
=  (f-g(e)h-eg^{2}(e),h+2eg(e),e)\circ (f,h,e) =V \circ  \Delta_g \circ V.
\ees
Hence 
\begin{equation}\label{7}
\Delta_g\circ S = S \circ \Delta_{-g}, \ \Delta_g \circ V=V\circ \tau \circ \Delta_{-g} \circ \tau.
\end{equation}

Consider the case $\lambda\neq 0$. In this case $G'$ is generated by $S$, $V$, and the automorphisms of the form 
\eqref{5}. Using \eqref{f6}, \eqref{6},  and \eqref{7}, we can represent any element of $G'$ in the form $\omega\circ\varphi$, where $\varphi\in G_1$ and $\omega$ is equal to $S$, $V$, or $id$. Notice that $\varphi\in G$, $V\circ\varphi\notin G$, and $S\circ\varphi\notin G$. It follows that $G=G_1$.  

Now consider the case $\lambda= 0$. We have
\bes
R_{\nu_1}\circ R_{\nu_2}=R_{\nu_1\nu_2},\\
S\circ R_\nu=(\nu^2e,-\nu h,f)=R_\nu\circ S,\\
V\circ R_\nu=(\nu^2f,\nu h,e)=R_\nu \circ V \circ H_{\nu^2},\\
\Delta_g \circ  R_\nu= R_\nu \circ \Delta_{\nu g}.
\ees

Using these relations we can represent any element of $G'$ in the form $R_{\nu}\circ \psi$, where $\psi$ belongs to the subgroup generated by  the automorphisms of the form \eqref{5}, $V$, and $S$. As proven above, $\psi$ can be represented in the form $\psi=\omega\circ\varphi$, where $\varphi\in G_1$ and $\omega$ is equal to $S$, $V$, or $\mathrm{id}$. Notice that $\varphi\in G$, $R_{\nu}\circ V\circ\varphi\notin G$, $R_{\nu}\circ S\circ\varphi\notin G$ for any $\nu$, and $R_{\nu} \circ \varphi\notin G$ for any $\nu\neq 1$.  It follows that $G=G_1$. $\Box$

\begin{corol} \label{c1} The group of linear automorphisms $A$ is generated by $\tau$ and  the linear automorphisms of the form  \eqref{5}.
\end{corol}

For any $a\in P_{\lambda}$ denote by $\mathrm{ad}(a) : P_{\lambda} \rightarrow P_{\lambda}$ the adjoint operator of $a$, i.e., the inner derivation of $P_{\lambda}$ such that 
$\mathrm{ad}(a)(x)=\{a,x\}$ for any $x\in P_{\lambda}$. If $D$ is a locally nilpotent derivation of $P_{\lambda}$, then
\bes
\exp(D)=\sum_{i=0}^{\infty} \frac{1}{i!}D^i
\ees
is an automorphism of $P_{\lambda}$ and is called an {\em exponential} automorphism. 

\begin{lemma} \label{l2} $(1)$ For any $g\in K[x]$ the derivation  $\mathrm{ad}(g(f))$ is locally nilpotent. 

$(2)$ Let $\Delta_g$ be an automorphism of the form \eqref{5} and let $\hat{g}=\int_0^x g dx$ be the formal integral of $g$. Then  
\bes
\Delta_g=\exp \mathrm{ad}(\hat{g}(f)).
\ees
\end{lemma}
\Proof 
It is easy to check that $ \{f^{n},h\}=2nf^{n},  \{f^{n},e\}=-nf^{n-1}h$. These relations easily imply the first statement of the lemma.

Suppose that $g(x)=\sum^{n}_{i=0} \alpha_{i} x^i$. Then $\hat{g}(x)=\sum^{n}_{i=0} \frac{1}{i+1}\alpha_{i} x^{i+1}$ and 
\bes
\hat{g}(f)=\sum^{n}_{i=0} \frac{1}{i+1}\alpha_{i} f^{i+1}.
\ees
Direct calculation gives that 
\bes
\{\hat{g}(f), h\}=2\sum^{n}_{i=0}\alpha_{i}f^{i+1}=2fg(f),\\
\{\hat{g}(f), e\}=-h\sum^{n}_{i=0}\alpha_{i}f^{i}=-hg(f),
\ees
and
\bes
\{\hat{g}(f), -hg(f)\}=-2fg^{2}(f).
\ees
These relations directly imply the statement of the lemma. $\Box$

This lemma clarifies the relations (\ref{f7}). Denote by $J$ the subgroup of all exponential automorphisms of the form \eqref{5} of  $G=\mathrm{Aut}\,P_{\lambda}$. 
Denote by  $T=\mathrm{Tr}(P_{\lambda})$ the subgroup of  $G=\mathrm{Aut}\,P_{\lambda}$ generated 
by the subgroup $J$ and by the subgroup of hyperbolic rotations $H$. The subgroup $T$ will be called the subgroup of {\em triangular} automorphisms of $P_{\lambda}$.

\begin{lemma} \label{l3}  $$T= J \rtimes H.$$
\end{lemma}
\Proof Direct calculation gives that 
\bes
H_\nu\circ \Delta_g\circ H_\nu^{-1}=\Delta_{q}, \  \text{where} \ q=\nu^{-1}g(\nu^{-1}f).
\ees
 This implies that  $T=JH$ and $J$ is a normal subgroup of $T$. 
Obviously, $J\cap H=\{\mathrm{id}\}$. $\Box$

\begin{corol} Any triangular automorphism of the algebra $P_{\lambda}$ can be written uniquely in the form 
\bes
(\alpha e-\alpha g(f)h-\alpha fg^{2}(f), h+2 fg(f),\alpha^{-1} f), \  g(x)\in K[x], \ \alpha\in K^\ast.
\ees
\end{corol}

\begin{lemma} \label{l4} The system of elements
\bes
B_0=\{ \Delta_q \ | \ q(x)\in xK[x] \}
\ees
is a system of representatives of the left cosets of  $C$ in $T$. 
\end{lemma}
\Proof Let $\psi\in T$. By Lemma  \ref{l3},  $\psi=\Delta_g\circ H_\nu$. Let  $g(f)=q(f)+\alpha$, where $q(x)\in xK[x]$ and $\alpha\in K$. By  \eqref{f7}, we get 
\bes
\psi=\Delta_g\circ H_\nu=\Delta_q\circ\Delta_\alpha\circ H_\nu.
\ees
Since $\Delta_\alpha\circ H_\nu\in C$, it follows that $\psi\in \Delta_q\circ C$. 

Suppose that $\Delta_{q_{1}}\circ C=\Delta_{q_{2}}\circ C$  . This means that
\bes
\Delta^{-1}_{q_{1}}\circ \Delta_{q_{2}}\in C.
\ees
By \eqref{f7}, we have
\bes
\Delta^{-1}_{q_{1}}\circ \Delta_{q_{2}}=\Delta{-q_{1}}\circ \Delta_{q_{2}}=\Delta_{-q_{1}+q_{2}}.
\ees
Obviously, $\Delta_{-q_{1}+q_{2}}\in C$ if and only if $q_{1}=q_{2}$. $\Box$

\section{Amalgamated free product structure of $G$}
\hspace*{\parindent}

We introduce a linear order on the set of basis words (\ref{f3}) of the algebra $P_{\lambda}$. Put $f>h>e$. Let $u$ and $v$ be elements of the form (\ref{f3}). We say that $u<v$ if $\deg(u)<\deg(v)$ or $\deg(u)=\deg(v)$ and $u<v$ with respect to the lexicographic order. Every nonzero element $g\in P_{\lambda}$ can be uniquely represented in the form
\bes
g=\alpha_{1}g_{1}+\alpha_{2}g_{2}+...+\alpha_{m}g_{m},
\ees
where $g_{i}$ are elements the form (\ref{f3}), $0\neq \alpha_{i}\in K$ for all $i$ and $g_{1}>g_{2}>...>g_{m}$.  Then $g_{1}$ is called the \emph{leading word } of the element $g$ and will be denoted by $\bar{g}$.

\begin{lemma} \label{l5} If $u, v, w$ are elements the form (\ref{f3}) and  $u<v$, then $\overline{uw}<\overline{vw}$.
\end{lemma}
\Proof  The statement of the lemma is obviously true if $\deg(u)<\deg(v)$. It is also true if $u, v, w$ are all elements of the form $f^{m}h^{n}$ or all elements of the form $h^{n}e^{r}$. Let  $\deg(u)=\deg(v)$. Consider the following cases:

1) Let $u=f^{m_{1}}h^{n_{1}}$, $v=f^{m_{2}}h^{n_{2}}$ and $w=h^{n}e^{r}$. We have $m_{1}<m_{2}$ since $u<v$. By (\ref{f2}), 
\bes
\overline{uw}=f^{m_{1}-r}h^{n_{1}+n+2r}, \ \ \overline{vw}=f^{m_{2}-r}h^{n_{2}+n+2r} \ \text{ if } \ r\leq m_{1};\\
\overline{uw}=h^{n_{1}+n+2m_{1}}e^{r-m_{1}}, \ \ \overline{vw}=f^{m_{2}-r}h^{n_{2}+n+2r} \ \text{ if } \  m_{1}<r<m_{2};\\
\overline{uw}=h^{n_{1}+n+2m_{1}}e^{r-m_{1}}, \ \ \overline{vw}=h^{n_{2}+n+2m_{2}}e^{r-m_{2}} \ \text{ if } \ m_{2}\leq r.
\ees 
Consequently, $\overline{uw}<\overline{vw}$.

2) Let  $u=h^{n_{1}}e^{r_{1}}$, $v=h^{n_{2}}e^{r_{2}}$ and $w=f^{m}h^{n}$. We have $r_{1}>r_{2}$ since $u<v$. By (\ref{f2}), 
\bes
\overline{uw}=h^{n_{1}+n+2m}e^{r_{1}-m}, \ \ \overline{vw}=h^{n_{2}+n+2m}e^{r_{2}-m} \ \text{ if } \ r_{2}\geq m; \\
\overline{uw}=h^{n_{1}+n+2m}e^{r_{1}-m}, \ \ \overline{vw}=f^{m-r_{2}}h^{n_{2}+n+2r_{2}} \ \text{ if } \ r_{1}>m>r_{2}; \\
\overline{uw}=f^{m-r_{1}}h^{n_{1}+n+2r_{1}}, \ \ \overline{vw}=f^{m-r_{2}}h^{n_{2}+n+2r_{2}} \ \text{ if } \ m\geq r_{1}. 
\ees 
Consequently, $\overline{uw}<\overline{vw}$.

3) Let  $u=f^{m_{1}}h^{n_{1}}$ and $v=h^{n_{2}}e^{r_{2}}$. This case contradicts the condition $u<v$.

4) Let $u=h^{n_{1}}e^{r_{1}}$, $v=f^{m_{2}}h^{n_{2}}$ and $w=h^{n}e^{r}$. We have $(r_{1},m_{2})\neq (0,0)$ since $u<v$. By (\ref{f2}),
\bes
\overline{uw}=h^{n_{1}+n}e^{r_{1}+r},\\
\overline{vw}=f^{m_{2}-r}h^{n_{2}+n+2r} \ \text{if} \ m_{2}\geq r,\\
\overline{vw}=h^{n_{2}+n+2m_{2}}e^{r-m_{2}} \ \text{if} \ m_{2}< r.
\ees
Consequently, $\overline{uw}<\overline{vw}$.

5) Let $u=h^{n_{1}}e^{r_{1}}$, $v=f^{m_{2}}h^{n_{2}}$ and $w=f^{m}h^{n}$. We have $(r_{1},m_{2})\neq (0,0)$ since $u<v$. By (\ref{f2}), 
\bes
\overline{uw}=f^{m-r_{1}}h^{n_{1}+n+2r_{1}} \ \text{if} \ m\geq r_{1},\\
\overline{uw}=h^{n_{2}+n+2m}e^{r_{1}-m} \ \text{if} \ m< r_{1},\\
\overline{vw}=f^{m_{2}+m}h^{n_{2}+n}.
\ees
Consequently, $\overline{uw}<\overline{vw}$. $\Box$

\begin{corol} \label{c2} Let $0\neq g,q\in P_{\lambda}$. Then $\overline{gq}=\overline{\bar{g} \bar{q}}$ and $\deg(gq)=\deg(g)+\deg(q)$.
\end{corol}

Let $G$ be an arbitrary group, $G_0, G_1, G_2$ be subgroups of the group $G$ and $G_0=G_1\cap G_2$. The group $G$ is called \textit{the free product of the subgroups $G_1$ and $G_2$ with the amalgamated subgroup $G_0$} and is denoted by $G=G_1\ast_{G_0} G_2$ if 
\begin{enumerate}
\item[(a)] $G$ is generated by the subgroups $G_1$ and $G_2$;
\item[(b)] the defining relations of the group $G$ consist only of the defining relations of the subgroups $G_1$ and $G_2$.
\end{enumerate}

If $S_1$ is a system of left representatives of $G_1$ by $G_0$ and $S_2$ is a system of left representatives of $G_2$ by $G_0$ then the group $G$ is a free product of subgroups $G_1$ and $G_2$ with the amalgamated subgroup $G_0$ (see, for example, \cite{Mag}) if and only if each  $g\in G$ is uniquely representable in the form 
$$g=g_1\ldots g_k c,$$
where $g_i\in S_1 \cup S_2,\: i=1,\ldots,k,$ $g_i, g_{i+1}$ do not belong to $S_1$ and $S_2$ at the same time, and $c\in G_0$.

For any automorphism $\phi=(f_{1},f_{2},f_{3})$ of the algebra $P_{\lambda}$ define its  multidegree by
\bes
\mathrm{mdeg}(\phi)= (\deg(f_{1}),\deg(f_{2}),\deg(f_{3})).
\ees
\begin{lemma} \label{l6} Let $A_0$ and  $B_0$ be the sets defined in Lemma \ref{l1} and Lemma \ref{l4}, respectively. Then any automorphism $\phi$ of the algebra $P_{\lambda}$  decomposes into a product of the form
\begin{gather} \label{9}
\phi=\alpha_1\circ \beta_1 \circ \alpha_2 \circ \beta_2 \circ\ldots \circ \alpha_k \circ\beta_k\circ \alpha_{k+1}\circ \lambda,
\end{gather}
where $\alpha_i\in A_0,\: \alpha_2,\ldots,\alpha_k\neq id$, $\beta_i\in B_0,\:\beta_1,\ldots,\beta_k\neq id$, $\lambda\in C$.
\end{lemma}
\Proof By Proposition \ref{p1},  any automorphism $\phi$ can be represented as
\bes
\phi= l_1 \circ t_1\circ l_2 \circ t_2\circ...\circ l_n \circ t_n \circ l_{n+1},
\ees
where  $l_{i}\in A$, $t_{i}\in T$, $t_{i}\notin A$.

We prove by induction on $n$ that $\phi$ can be represented as \eqref{9} with $k\leq n$. By Lemma \ref{l1}, we have  $l_1=\alpha_{1}\circ \lambda_{1}$, where $\alpha_1 \in A_0, \lambda_1\in C$. By \eqref{f7}, $\lambda_1\circ t_{1} \in T$ and $\lambda_1\circ t_{1} \notin A$. Then, by Lemma \ref{l4}, we have $\lambda_1\circ t_{1}=\beta_1\circ\lambda'_1$, where $id\neq\beta_1\in B_0$, $\lambda'_1\in C$. Hence
\bes
l_1 \circ t_1=\alpha_{1}\circ \beta_1\circ\lambda'_1.
\ees
Consequently, 
\bes
\phi= \alpha_{1}\circ \beta_1\circ (\lambda'_1 \circ l_2) \circ t_2\circ...\circ l_n \circ t_n \circ l_{n+1}.
\ees
By the induction proposition, the product
$$
(\lambda'_1 \circ l_2) \circ t_2\circ...\circ l_n \circ t_n \circ l_{n+1}
$$
can be written in the form 
\bes
\alpha_2 \circ \beta_2 \circ\ldots \circ \alpha_k \circ\beta_k\circ \alpha_{k+1}\circ \lambda, \ \ k\leq n.
\ees
Consequently, 
\bes
\phi=\alpha_1\circ \beta_1 \circ \alpha_2 \circ \beta_2 \circ\ldots \circ \alpha_k \circ\beta_k\circ \alpha_{k+1}\circ \lambda.
\ees

If $\alpha_2 \neq id$, then the resulting representation has the form \eqref{9}. Let $\alpha_2 = id$. Since $\beta_1 \circ \beta_2=\beta'_2 \in B_0$, it follows that
\bes
\phi=\alpha_1\circ \beta_1 \circ \beta_2 \circ\ldots \circ \alpha_k \circ\beta_k\circ \alpha_{k+1}\circ \lambda=\alpha_1\circ \beta'_2 \circ\ldots \circ \alpha_k \circ\beta_k\circ \alpha_{k+1}\circ \lambda.
\ees
Consequently, $\phi$ can be represented as \eqref{9} since $k-1<n$. $\Box$

\begin{lemma} \label{l7} Let $\varphi=(f_1,f_2,f_3)$  be an automorphism of the algebra $P_{\lambda}$ such that  
\bes
\varphi=(f_1,f_2,f_3)=\beta_1\circ\alpha_2\circ\beta_2\circ\ldots\circ\alpha_k\circ\beta_k,
\ees
where $id\neq \alpha_i\in A_0$ and $id\neq \beta_i\in B_0$ for any $i$.
If $\beta_i=\Delta_{q_{i}}$, where $q_i(x)\in xK[x]$  and $\deg(q_i)=n_i>0$, then
\bes
\deg(f_1)=\prod^{k}_{i=1}(2n_{i}+1), \ \ 
\deg(f_2)=(n_{k}+1)\prod^{k-1}_{i=1}(2n_{i}+1), \ \ 
\deg(f_3)= \prod^{k-1}_{i=1}(2n_{i}+1).
\ees
(here we assume  that $\prod^{k-1}_{i=1}(2n_{i}+1)=1$   if $k = 1$).
\end{lemma}
\Proof Prove the statement of the lemma by induction on $k$. If $k=1$, then $\varphi=\beta_1$ and $\mathrm{mdeg}(\varphi)=(2n_1+1,n_1+1,1)$.

Suppose that the lemma holds for $k-1$. Put
\bes
\varphi_1=\beta_1\circ\alpha_2\circ\beta_2\circ\ldots\circ\alpha_{k-1}\circ\beta_{k-1}=(g_1,g_2,g_3).
\ees
By  the induction  hypothesis, we get
\bes
\deg(g_1)=\prod^{k-1}_{i=1}(2n_{i}+1), \ \ \ 
\deg(g_2)=(n_{k-1}+1)\prod^{k-2}_{i=1}(2n_{i}+1), \ \ \ 
\deg(g_3)=\prod^{k-2}_{i=1}(2n_{i}+1).
\ees
Then
\bes
\varphi=(f_1,f_2,f_3)=\beta_1\circ\alpha_2\circ\beta_2\circ\ldots\circ\alpha_k\circ\beta_k=\varphi_1\circ\alpha_k\circ\beta_k=(g_1,g_2,g_3)\circ\alpha_k\circ\beta_k.
\ees
We have  
\bes
(u_1,u_2,u_3)=(g_1,g_2,g_3)\circ\alpha_k=(g_3, -g_2+2\beta g_3, g_1+\beta g_2-\beta^{2} g_{3}).
\ees
Since $\deg (g_{1})>\deg(g_{2})>\deg(g_{3})$, it follows that
\bes
\deg(u_1)=\deg(g_3)=\prod^{k-2}_{i=1}(2n_{i}+1),\\
\deg(u_2)=\deg(g_2)=(n_{k-1}+1)\prod^{k-2}_{i=1}(2n_{i}+1),\\
\deg(u_3)=\deg(g_1)=\prod^{k-1}_{i=1}(2n_{i}+1).
\ees
Further,
\bes
\varphi=(f_1,f_2,f_3)=(u_1,u_2,u_3)\circ\beta_k=(u_1,u_2,u_3)\circ (e-q_{k}(f)h-fq^{2}_{k}(f), h+2fq_{k}(f), f)=\\
=(u_1-q_{k}(u_3)u_2-u_3q^{2}_{k}(u_3), u_2+2u_3q_{k}(u_3), u_3).
\ees
Consequently,
\bes
\deg(f_1)=max\{\deg(u_1),\deg(q_k(u_3)u_2),\deg(u_3q^{2}_k(u_3))\},\\
\deg(f_2)=max\{\deg(u_2),\deg(u_3 q_k(u_3))\},\\
\deg(f_3)=\deg(u_3).
\ees
Recall that  $\deg(q_k)=n_k>0$ and $\deg(u_3)=\prod^{k-1}_{i=1}(2n_{i}+1).$
Then, by Corollary \ref{c2}, 
\bes
\deg(u_3 q_k(u_3))=\prod^{k-1}_{i=1}(2n_{i}+1)+n_{k}\prod^{k-1}_{i=1}(2n_{i}+1)=(n_{k}+1)\prod^{k-1}_{i=1}(2n_{i}+1),\\
\deg(q_k(u_3)u_2)=n_k\prod^{k-1}_{i=1}(2n_{i}+1)+(n_{k-1}+1)\prod^{k-2}_{i=1}(2n_{i}+1),\\
\deg(u_3q^{2}_k(u_3))= \prod^{k-1}_{i=1}(2n_{i}+1)+2n_k\prod^{k-1}_{i=1}(2n_{i}+1)=\prod^{k}_{i=1}(2n_{i}+1).
\ees
Consequently,
\bes
\deg(f_1)=\prod^{k}_{i=1}(2n_{i}+1), \ \ 
\deg(f_2)=(n_{k}+1)\prod^{k-1}_{i=1}(2n_{i}+1), \ \ 
\deg(f_3)= \prod^{k-1}_{i=1}(2n_{i}+1). \ \Box
\ees

\begin{lemma} \label{l8} The decomposition \eqref{9} of the automorphism $\phi$ from Lemma \ref{l6} is unique.
\end{lemma}
\Proof It suffices to show that
\bes
\alpha_1\circ\beta_1\circ\alpha_2\circ\beta_2\circ\ldots\circ\alpha_k\circ \beta_k\circ\alpha_{k+1}\circ\lambda\neq id,
\ees
for $k\geq 1$, $\alpha_i\in A_0,\:\alpha_2,\ldots,\alpha_k\neq id$, $\beta_i\in B_0,\:\beta_1,\ldots,\beta_k\neq id$, $\lambda\in C$.

Suppose that
\bes
\alpha_1\circ\beta_1\circ\alpha_2\circ\beta_2\circ\ldots\circ\alpha_k\circ\beta_k\circ\alpha_{k+1}\circ\lambda=\mathrm{id}.
\ees
Then
\begin{gather} \label{10}
\beta_1\circ\alpha_2\circ\beta_2\circ\ldots\circ\alpha_k\circ\beta_k=\alpha_1^{-1}\circ\lambda^{-1}\circ\alpha_{k+1}^{-1}.
\end{gather}
By Lemma \ref{l7}, the automorphism 
\bes
\varphi=\beta_1\circ\alpha_2\circ\beta_2\circ\ldots\circ\alpha_k\circ\beta_k
\ees
has the multidegree
\bes
\mathrm{mdeg}(\varphi)=\left( \prod^{k}_{i=1}(2n_{i}+1), (n_{k}+1)\prod^{k-1}_{i=1}(2n_{i}+1), \prod^{k-1}_{i=1}(2n_{i}+1)\right),
\ees
where $n_{i}>0$. The automorphism on the right-hand side of \eqref{10} is linear and has the multidegree $(1,1,1)$.
This contradicts \eqref{10}. $\Box$

\begin{theorem} \label{th2}  The group $G$ of automorphism of the Poisson algebra $P_{\lambda}$ is the free product of subgroups $A$ and $T$ with amalgamated subgroup $C=A \cap T$, i.e.,
\bes
G = A\ast_ {C} T.
\ees 
\end{theorem}
\Proof Recall that $A_0$ and $B_0$ are  the systems of representatives of the left cosets of  $C$ in $A$ and $T$, respectively. By Lemma \ref{l6} and Lemma \ref{l8}, any automorphism of $P_{\lambda}$ can be uniquely represented as \eqref{9}. According to \cite{Mag}, we have 
$$ G = A\ast_ {C} T. \ \ \Box $$

\section{Isomorphism of automorphism groups}
\hspace*{\parindent}

For any $\lambda\in K$ let  $U_{\lambda}=U(\mathrm{sl}_{2}(K))/(C_U-\lambda)$ be the quotient algebra of the  algebra $U(\mathrm{sl}_{2}(K)$ by the principal ideal $(C_U-\lambda)$. Denote by  $e,h,f$  the images of $E,H,F$ in $U_{\lambda}$, respectively. Then we have
$$\left[e,f\right]=h, \left[f,h\right]=2f, \left[h,e\right]=2e,\lambda=4fe+h^{2}+2h=4ef+h^{2}-2h=2ef+2fe+h^{2}.$$

J. Dixmier \cite{Dixmier73} described generators of the group of automorphisms of the algebra $U_{\lambda}$.

\begin{theorem} \label{th3} \cite[p. 563]{Dixmier73} The group $\mathrm{Aut}(U_{\lambda})$ is generated by the exponential automorphisms $\Phi_{n,\mu}$ and $\Psi_{n,\mu}$ ($n>0, \mu\in K$), where 
\bes
\Phi_{n,\mu}=(e-\mu nf^{n-1}h+\mu n(n-1)f^{n-1}-\mu^{2}n^{2}f^{2n-1}, h+2\mu nf^{n},  f),\\
\Psi_{n,\mu}=(e, h-2\mu ne^{n}, f+\mu ne^{n-1}h+\mu n(n-1)e^{n-1}-\mu^{2}n^{2}e^{2n-1}).
\ees 
\end{theorem}

Recall that $A=\mathrm{Aut}(\mathrm{sl}_{2}(K))$ is the group of all automorphisms of the Lie algebra $\mathrm{sl}_{2}(K)$. In Section 2 we identified $A$ with a subgroup of 
 automorphisms  $G=\mathrm{Aut}(P_{\lambda})$ of the Poisson algebra $P_{\lambda}$. Obviously, every automorphism of $A$ gives a unique automorphism of the algebra 
$U_{\lambda}$. For this reason we identify $A$ with the corresponding subgroup of the group of automorphisms of $U_{\lambda}$. Thus, without loss of generality, we can assume that $\mathrm{Aut}(P_{\lambda})$ and $\mathrm{Aut}(U_{\lambda})$ both contain $A$. Then the subgroups $H$ and $C$ of $A$ may be considered as subgroups of $\mathrm{Aut}(P_{\lambda})$ and $\mathrm{Aut}(U_{\lambda})$.

Denote by $J'$ the subgroup of all exponential automorphisms $\Phi_{n,\mu}$ ($\nu \in K^{\ast},  n\in \mathbb{N}^{\ast},\mu \in K$) of  $\mathrm{Aut}\,U_{\lambda}$. 
Notice that every automorphism $\Phi_{n,\mu}$ can be written in the form 
\begin{equation}\label{12}
\begin{split}
\delta_g(e)=e-g(f)h-fg^{2}(f)+fg'(f),\\
\delta_g(h)=h+2fg(f),\\
\delta_g(f)=f, 
\end{split}
\end{equation} 
where $g(x)\in K[x]$ and $g'(x)$ is the formal derivative of $g(x)$. 

It is easy to check that
\bee\label{13}
\delta_{g_1}\circ \delta_{g_2}=\delta_{g_1+g_2}. 
\eee

Denote by  $T'=\mathrm{Tr}(P_{\lambda})$ the subgroup of  $\mathrm{Aut}\,U_{\lambda}$ generated 
by the subgroup $J'$ and by the subgroup of hyperbolic rotations $H$. The subgroup $T'$ will be called the subgroup of {\em triangular} automorphisms of $U_{\lambda}$.

O. Fleury  \cite{Fleury}  proved that the automorphism group of the algebra $U_{\lambda}$  admits an amalgamated free product structure.

\begin{theorem} \label{th4} \ \cite{Fleury} We have
\bes
\mathrm{Aut}(U_{\lambda})=A\ast_{C} T',
\ees
where $A=\mathrm{Aut}(\mathrm{sl}_2(K))$ is the subgroup of all automorphisms of the Lie algebra $\mathrm{sl}_{2}(K)$ in $\mathrm{Aut}(U_{\lambda})$ and $ C=A\cap T'$.
\end{theorem}

\begin{lemma} \label{l9}  $$T'= J' \rtimes H.$$
\end{lemma}
\Proof Direct calculation gives that 
\bes
H_\nu\circ \delta_g\circ H_\nu^{-1}=\delta_{q}, \ \ \text{where} \ q=\nu^{-1}g(\nu^{-1}f).
\ees
This gives that  $T'=J'H$ and $J'$ is a normal subgroup of $T'$. 
Obviously, $J'\cap H=\{\mathrm{id}\}$. $\Box$

\begin{corol} Any triangular automorphism of the algebra $U_{\lambda}$ can be written uniquely in the form 
\bes
(\alpha e-\alpha g(f)h-\alpha fg^{2}(f)+\alpha fg'(f), h+2 fg( f),\alpha^{-1} f), \  g(x)\in K[x], \ \alpha\in K^\ast.
\ees
\end{corol}

\begin{theorem} \label{th5} The automorphism groups of the algebras $P_{\lambda}$ and  $U_{\lambda}$ are isomorphic, i.e., 
\bes
\mathrm{Aut}(P_{\lambda})\cong \mathrm{Aut}(U_{\lambda}).
\ees
\end{theorem}
\Proof Consider the map $\phi: J'\rightarrow J $  given by the rule  $\delta_g \mapsto \Delta_g$. By \eqref{f7} and \eqref{13},
\bes
\phi(\delta_{g_1}\circ\delta_{g_2})=\phi(\delta_{g_1})\circ \phi(\delta_{g_2}).
\ees
Consequently, $J'\cong J$. By Lemma \ref{l3} and Lemma \ref{l9}, we get $T'\cong T$. Then 
Theorem \ref{th2} and Theorem \ref{th4}  directly imply the statement of the theorem. $\Box$

\bigskip

{\bf Acknowledgments}

\bigskip
The first author was supported by the grant of the Ministry of Education and Science of the Republic of Kazakhstan (project AP08052290).

The second author was supported by the grant of the Ministry of Education and Science of the Republic of Kazakhstan (project  AP09261086).

\end{document}